\documentclass[10pt]{article}
\usepackage{cite}
\usepackage{mathrsfs}
\usepackage{amsfonts}
\usepackage{amsmath}
\usepackage{amsfonts,amssymb,color}
\usepackage{dsfont}
\usepackage{curves}
\usepackage{mathrsfs}
\usepackage{pifont}
\usepackage{amssymb}
\allowdisplaybreaks

\numberwithin{equation}{section}

\date{}

\def\BigRoman{\uppercase\expandafter{\romannumeral\number\count 255 }}
\def\Romannumeral{\afterassignment\BigRoman\count255=}

\begin{document}
\title{The $A_{\alpha}$-spectral radius for a matching-$k$-cover in a graph
}
\author{\small Quanru Pan, Sizhong Zhou\footnote{Corresponding author. E-mail address: zsz\_cumt@163.com}\\
\small School of Science, Jiangsu University of Science and Technology,\\
\small Zhenjiang, Jiangsu 212100, China\\
}

\maketitle
\begin{abstract}
\noindent Let $G$ be a graph with vertex set $V(G)$ and edge set $E(G)$. For $\alpha\in[0,1]$, Nikiforov [Appl. Anal. Discrete Math. 11 (2017) 81--107] defined the $A_{\alpha}$-matrix of $G$ as
$A_{\alpha}(G)=\alpha D(G)+(1-\alpha)A(G)$, where $D(G)$ and $A(G)$ respectively denote the degree diagonal matrix and the adjacency matrix of $G$. The largest eigenvalue of $A_{\alpha}(G)$ is
said to be the $A_{\alpha}$-spectral radius of $G$, and written as $\rho_{\alpha}(G)$. A set of independent edges in $G$ is called a matching. A matching-$k$-cover in $G$ is a collection of $k$
matchings whose union covers $V(G)$. Ait Ferhat et al. [Math. Program. 203 (2024) 271--284] proved that a connected graph $G$ contains a matching-$k$-cover if and only if $|S|\leq k|N_G(S)|$ for
each independent subset $S$ of $V(G)$, where $k\geq2$ is an integer. Utilizing the above result, Shang, Zhang and Wang [Discrete Appl. Math. 398 (2027) 498--505] obtained a sufficient condition
based on the adjacency spectral radius for a connected graph to possess a matching-$k$-cover for $k\geq2$. In this paper, we propose a tight sufficient condition in terms of the $A_{\alpha}$-spectral
radius to guarantee the existence of a matching-$k$-cover in a connected graph for $k\geq2$, which generalizes Shang, Zhang and Wang's result. As a corollary, we also present a signless Laplacian spectral
radius condition for a connected graph to contain a matching-$k$-cover for $k\geq2$.
\\
\begin{flushleft}

(2020) Mathematics Subject Classification: 05C50, 05C70
\end{flushleft}
\end{abstract}

\section{Introduction}

In this paper, we only deal with finite undirected simple graphs. Let $G=(V(G),E(G))$ be a graph, where $V(G)$ and $E(G)$ respectively denote the vertex set and the edge set of $G$. The order of $G$
is defined as $n=|V(G)|$. For $v\in V(G)$, we denote by $N_G(v)$ the neighborhood of $v$ in $G$, and by $d_G(v)$ the degree of $v$ in $G$. A vertex $v$ in $G$ is called a pendant vertex if $d_G(v)=1$.
We write $N_G(S)=\bigcup\limits_{v\in S}N_G(v)$ for any $S\subseteq V(G)$. For $S\subseteq V(G)$, we denote by $G[S]$ the subgraph of $G$ induced by $S$, and by $G-S$ the subgraph of $G$ induced by
$V(G)\setminus S$. A vertex subset $X$ in $G$ is called independent if $G[X]$ has no edges. A complete graph of order $n$ is denoted by $K_n$. Let $G_1$ and $G_2$ be two vertex-disjoint graphs. The
the union of $G_1$ and $G_2$, denoted by $G_1\cup G_2$, is the graph with vertex set $V(G_1)\cup V(G_2)$ and edge set $E(G_1)\cup E(G_2)$. The join of $G_1$ and $G_2$, denoted by $G_1\vee G_2$, is
the graph obtained from $G_1\cup G_2$ by adding all edges joining vertices of $G_1$ to vertices of $G_2$.

A set of independent edges in $G$ is called a matching. A perfect matching in $G$ is a matching such that each vertex of $G$ is incident with precisely one edge in it. A matching-$k$-cover in $G$ is
a collection of $k$ matchings whose union covers $V(G)$. Obviously, $G$ has a matching-$k'$-cover if it has a matching-$k$-cover for every integer $k'\geq k$. In fact, a matching-1-cover is precisely
a perfect matching.

For a given graph $G$ with vertex set $V(G)=\{v_1,v_2,\ldots,v_n\}$, the adjacency matrix of $G$, denoted by $A(G)$, is defined as $A(G)=(a_{ij})_{n\times n}$, where $a_{ij}=1$ if $v_iv_j\in E(G)$,
and $a_{ij}=0$ otherwise. Write $D(G)=\mbox{diag}(d_G(v_1),d_G(v_2),\ldots,d_G(v_n))$. The signless Laplacian matrix of $G$, denoted by $Q(G)$, is defined by $Q(G)=D(G)+A(G)$. For $\alpha\in[0,1]$,
Nikiforov \cite{N} defined the $A_{\alpha}$-matrix of $G$ as $A_{\alpha}(G)=\alpha D(G)+(1-\alpha)A(G)$. In fact, $A_0(G)=A(G)$ and $A_{\frac{1}{2}}(G)=\frac{1}{2}Q(G)$. The eigenvalues of $A_{\alpha}(G)$
are said to be the $A_{\alpha}$-eigenvalues of $G$, and the largest one of them is called the $A_{\alpha}$-spectral radius of $G$, and written as $\rho_{\alpha}(G)$. The adjacency and signless Laplacian
eigenvalues (and spectral radius) can be similarly defined. We denote by $\rho(G)$ the adjacency spectral radius of $G$, and by $q(G)$ the signless Laplacian spectral radius of $G$. We refer the reader
to \cite{PR,So,LL,JLZ,ZLB,Zs1,Zs,Za,WZL,Ws,Wu} for some properties of the $A_{\alpha}$-spectral radius.

Hall \cite{H} verified that a bipartite graph with bipartition $(X,Y)$ contains a matching saturating $X$ if and only if $|N_G(S)|\geq|S|$ occurs for each $S\subseteq X$. In particular, it gives a
characterization for a bipartite graph with a perfect matching. For a general graph, Tutte \cite{T} provided a necessary and sufficient condition for a graph with a perfect matching. Based on Tutte's
result, the relationships between perfect matchings and the spectral radius of various graph matrices have been studied by some scholars. In 2021, O \cite{Os} posed a sufficient condition based on the
adjacency spectral radius for a connected graph with a perfect matching. Fan and Lin \cite{FL}, Wu and Zhou \cite{WZ} proposed an adjacency spectral radius condition for the existence of a perfect matching in a 1-binding graph.
Liu, Fan and Shu \cite{LFS} put forward two tight sufficient conditions according to the adjacency spectral radius and the signless Laplacian spectral radius for a tough graph to possess a perfect matching.
Zhang and Lin \cite{ZL1} provided two upper bounds for the distance spectral radius which ensure the existence of a perfect matching in a graph or a bipartite graph. Hu and Zhang \cite{HZ} created a
connection between perfect matchings and the distance spectral radii in 1-binding graphs. Zhao, Huang and Wang \cite{ZHW} provided an $A_{\alpha}$-spectral radius condition for a graph to possess a perfect
matching. Zhou and Lin \cite{ZL2} showed a sufficient condition in view of the $A_{\alpha}$-spectral radius for a 1-binding graph to admit a perfect matching.

Note that a matching-$k$-cover is a natural generalization of a perfect matching: instead of utilizing one matching to cover all vertices of a graph, we allow the vertices to be covered by the union
of $k$ matchings. Wang, Song and Yuan \cite{WSY} studied the optimal matching covers in graphs. Ait Ferhat et al. \cite{AKSS} showed a characterization for a graph with a matching-$k$-cover for $k\geq2$.
Shang, Zhang and Wang \cite{SZW} obtained a sufficient condition based on the adjacency spectral radius for a connected graph to possess a matching-$k$-cover for $k\geq2$.

\medskip

\noindent{\textbf{Theorem 1.1}} (Shang, Zhang and Wang \cite{SZW}) Let $k\geq2$ be an integer, and let $G$ be a connected graph of order $n$ with $n\geq k+2$.

(\romannumeral1) For $k=2$ and $n=7$, if $\rho(G)\geq\rho(K_2\vee5K_1)$, then $G$ admits a matching-2-cover, unless $G=K_2\vee5K_1$.

(\romannumeral2) For $k=3$ and $n=9$, if $\rho(G)\geq\rho(K_2\vee7K_1)$, then $G$ admits a matching-3-cover, unless $G=K_2\vee7K_1$.

(\romannumeral3) For $k=2$ and $n\neq7$, or $k=3$ and $n\neq9$, or $k\geq4$, if $\rho(G)\geq\rho(K_1\vee((k+1)K_1\cup K_{n-k-2}))$, then $G$ admits a matching-$k$-cover, unless $G=K_1\vee((k+1)K_1\cup K_{n-k-2})$.

\medskip

Motivated by Theorem 1.1, we propose a tight sufficient condition in terms of the $A_{\alpha}$-spectral radius to guarantee the existence of a matching-$k$-cover in a connected graph, which generalizes
Theorem 1.1 on $\alpha=0$ due to Shang, Zhang and Wang \cite{SZW}.

\medskip

\noindent{\textbf{Theorem 1.2.}} Let $\alpha\in[0,\frac{3}{4}]$, and $k\geq2$ be an integer. Let $G$ be a connected graph of order $n$ with $n\geq 4k+12$. If
$$
\rho_{\alpha}(G)\geq\rho_{\alpha}(K_1\vee(K_{n-k-2}\cup(k+1)K_1)),
$$
then $G$ contains a matching-$k$-cover, unless $G=K_1\vee(K_{n-k-2}\cup(k+1)K_1)$.

\medskip

The following corollary occurs by setting $\alpha=\frac{1}{2}$ in Theorem 1.2.

\medskip

\noindent{\textbf{Corollary 1.3.}} Let $k\geq2$ be an integer, and let $G$ be a connected graph of order $n$ with $n\geq 4k+12$. If
$$
q(G)\geq q(K_1\vee(K_{n-k-2}\cup(k+1)K_1)),
$$
then $G$ contains a matching-$k$-cover, unless $G=K_1\vee(K_{n-k-2}\cup(k+1)K_1)$.

\section{Some preliminaries}

In this section, we list some essential lemmas which will be used to verify our main results. Ait Ferhat et al. \cite{AKSS} arose a characterization for a graph containing a matching-$k$-cover.

\medskip

\noindent{\textbf{Lemma 2.1}} (Ait Ferhat et al. \cite{AKSS}). Let $G$ be a connected graph. Then $G$ contains a matching-$k$-cover if and only if
$$
|S|\leq k|N_G(S)|
$$
for each independent subset $S$ of $V(G)$, where $k\geq2$ is an integer.

\medskip

\noindent{\textbf{Lemma 2.2}} (Nikiforov \cite{N}). Let $H$ be a subgraph of a connected graph $G$. Then
$$
\rho_{\alpha}(G)\geq\rho_{\alpha}(H),
$$
with equality following if and only if $G=H$.

\medskip

Let $M$ be a real $n\times n$ matrix, and let $\pi:\mathcal{N}=\mathcal{N}_1\cup\mathcal{N}_2\cup\cdots\cup\mathcal{N}_r$ be a partition of $\mathcal{N}=\{1,2,\ldots,n\}$. Based on this partition $\pi$, the matrix
$M$ can be correspondingly written in block form as
\begin{align*}
M=\left(
  \begin{array}{cccc}
    M_{11} & M_{12} & \cdots & M_{1r}\\
    M_{21} & M_{22} & \cdots & M_{2r}\\
    \vdots & \vdots & \ddots & \vdots\\
    M_{r1} & M_{r2} & \cdots & M_{rr}\\
  \end{array}
\right),
\end{align*}
where $M_{ij}$ denotes the submatrix of $M$ with rows indexed by $\mathcal{N}_i$ and columns indexed by $\mathcal{N}_j$. The matrix $M_{\pi}=(m_{ij})_{r\times r}$ is called the quotient matrix of $M$ based on $\pi$,
where $m_{ij}$ is the average row sum of $M_{ij}$. The partition $\pi$ is called equitable If each row of every block $M_{ij}$ admits the same sum.

\medskip

\noindent{\textbf{Lemma 2.3}} (You et al. \cite{YYSX}). Let $M$ be a real $n\times n$ matrix with an equitable partition $\pi$, and let $M_{\pi}$ be the corresponding quotient matrix. Then the eigenvalues
of $M_{\pi}$ are also eigenvalues of $M$. Furthermore, if $M$ is nonnegative and irreducible, then the largest eigenvalues of $M$ and $M_{\pi}$ are equal.

\medskip

A matrix $M$ is said to be a Hermitian matrix if $M=(\overline{M})^{T}$, where $(\overline{M})^{T}$ denotes the conjugate transpose of $M$. It is obvious that a real symmetric matrix is a Hermitian matrix. Next, we
show the well-known Cauchy Interlacing Theorem.

\medskip

\noindent{\textbf{Lemma 2.4}} (Haemers \cite{Hi}). Let $M$ be a Hermitian matrix of order $s$, and let $N$ be a principal submatrix of $M$ with order $t$. If $\lambda_1\geq\lambda_2\geq\cdots\geq\lambda_s$ are the
eigenvalues of $M$ and $\mu_1\geq\mu_2\geq\cdots\geq\mu_t$ are the eigenvalues of $N$, then $\lambda_i\geq\mu_i\geq\lambda_{s-t+i}$ for $1\leq i\leq t$.

\section{The proof of Theorem 1.2}

\noindent{\it Proof of Theorem 1.2.} Suppose that $G$ contains no matching-$k$-cover. Among all such graphs, select $G$ with maximum $A_{\alpha}$-spectral radius. Using Lemma 2.1, there exists an independent set
$S\subseteq V(G)$ satisfying $|S|\geq k|N_G(S)|+1$. Write $s=|S|$ and $p=|N_G(S)|$. Then $s\geq kp+1$.

Let $G_1$ be derived from $G$ by adding all possible edges while keeping $S$ independent and keeping $N_{G_1}(S)$ as the neighborhood of $S$. Then we conclude $G_1=K_p\vee(K_{n-s-p}\cup sK_1)$. Note that $S$ is
still an independent set of $G_1$ and $|S|\geq k|N_{G_1}(S)|+1$. Then Lemma 2.1 yields that $G_1$ contains no matching-$k$-cover. According to Lemma 2.2 and $G$ being a spanning subgraph of $G_1$, we infer
$\rho_{\alpha}(G)\leq\rho_{\alpha}(G_1)$. From the choice of $G$, we possess $\rho_{\alpha}(G)=\rho_{\alpha}(G_1)$. Notice that $G_1$ is connected. If $G$ is a proper subgraph of $G_1$, then it follows from Lemma
2.2 that $\rho_{\alpha}(G)<\rho_{\alpha}(G_1)$, a contradiction. Thus, we obtain
$$
G=G_1=K_p\vee(K_{n-s-p}\cup sK_1).
$$

If $s\geq kp+2$, select a subset $S_1\subseteq S$ with $|S_1|=kp+1$, and set $S_2=S\setminus S_1$. We construct a graph $G_2$ from $G_1$ by adding edges in $G_1-S_1$ so that $G_2-S_1$ is a complete graph. Then
$$
G_2=K_p\vee(K_{n-(k+1)p-1}\cup(kp+1)K_1).
$$
Note that $S_1$ is still an independent set of $G_2$ and $|N_{G_2}(S_1)|=p$. Thus, we get $|S_1|=kp+1>k|N_{G_2}(S_1)|$. Together with Lemma 2.1, $G_2$ contains no matching-$k$-cover. Moreover, $G_2$ contains $G_1$
as a proper spanning subgraph. Combining this with Lemma 2.2 and $G_2$ being connected, we obtain $\rho_{\alpha}(G_1)<\rho_{\alpha}(G_2)$. Recall that $G=G_1$. Then $\rho_{\alpha}(G)<\rho_{\alpha}(G_2)$ holds,
which contradicts the choice of $G$ as a counterexample with maximum $A_{\alpha}$-spectral radius. Hence, $s=kp+1$, and so
$$
G=K_p\vee(K_{n-(k+1)p-1}\cup(kp+1)K_1),
$$
where $n\geq(k+1)p+1$.

If $p=1$, then $G=K_1\vee(K_{n-k-2}\cup(k+1)K_1)$ and $\rho_{\alpha}(G)=\rho_{\alpha}(K_1\vee(K_{n-k-2}\cup(k+1)K_1)$. Observe that $K_1\vee(K_{n-k-2}\cup(k+1)K_1)$ contains no matching-$k$-cover, a contradiction.
Next, we assume that $p\geq2$.

Let $G_*=K_1\vee(K_{n-k-2}\cup(k+1)K_1)$. Then the quotient matrix of $A_{\alpha}(G_*)$ in view of the partition $V(G_*)=V(K_1)\cup V(K_{n-k-2})\cup V((k+1)K_1)$ is
\begin{align*}
B_*=\left(
  \begin{array}{ccc}
  \alpha n-\alpha & (1-\alpha)(n-k-2) & (1-\alpha)(k+1)\\
  1-\alpha & n+\alpha-k-3 & 0\\
  1-\alpha & 0 & \alpha\\
  \end{array}
\right).
\end{align*}
The characteristic polynomial of $B_*$ is equal to
\begin{align*}
\varphi_{B_*}(x)=&x^{3}+(k-\alpha+3-(1+\alpha)n)x^{2}+(\alpha n^{2}+(\alpha^{2}-k\alpha-\alpha-1)n-2\alpha+1)x\\
&-\alpha^{2}n^{2}+(2k\alpha^{2}-2k\alpha+k+3\alpha^{2}-\alpha+1)n-(k^{2}+5k+4)\alpha^{2}\\
&+(2k^{2}+8k+5)\alpha-k^{2}-4k-3.
\end{align*}
Using Lemma 2.3 and the equitable partition $V(G_*)=V(K_1)\cup V(K_{n-k-2})\cup V((k+1)K_1)$, the largest root of $\varphi_{B_*}(x)=0$ is equal to $\rho_{\alpha}(G_*)$. Write $\theta=\rho_{\alpha}(G_*)$. Then
$\varphi_{B_*}(\theta)=0$.

Notice that $G_*=K_1\vee(K_{n-k-2}\cup(k+1)K_1)$ contains $K_{n-k-1}$ as a proper subgraph. Combining this with Lemma 2.2, we conclude
\begin{align}\label{eq:3.1}
\theta=\rho_{\alpha}(G_*)>\rho_{\alpha}(K_{n-k-1})=n-k-2.
\end{align}
The following proof will be divided into three cases.

\noindent{\bf Case 1.} $n\geq(k+1)p+3$.

Recall that $G=K_p\vee(K_{n-(k+1)p-1}\cup(kp+1)K_1)$. Based on the partition $V(G)=V(K_p)\cup V(K_{n-(k+1)p-1})\cup V((kp+1)K_1)$, the quotient matrix of $A_{\alpha}(G)$ can be correspondingly written as
\begin{align*}
B=\left(
  \begin{array}{ccc}
  \alpha n+(1-\alpha)p-1 & (1-\alpha)(n-(k+1)p-1) & (1-\alpha)(kp+1)\\
  (1-\alpha)p & n+(\alpha-k-1)p-2 & 0\\
  (1-\alpha)p & 0 & \alpha p\\
  \end{array}
\right).
\end{align*}
Its characteristic polynomial is given by
\begin{align}\label{eq:3.2}
\varphi_B(x)=&x^{3}+((k-\alpha)p+3-(1+\alpha)n)x^{2}\nonumber\\
&+(\alpha n^{2}+(\alpha^{2}-k\alpha+\alpha)pn-(1+2\alpha)n-kp^{2}+(k-1-2\alpha)p+2)x\nonumber\\
&-\alpha^{2}pn^{2}+(2k\alpha^{2}-2k\alpha+k)p^{2}n+(3\alpha^{2}-\alpha+1)pn\nonumber\\
&-((k^{2}+2k)\alpha^{2}-(2k^{2}+3k)\alpha+k^{2}+k)p^{3}\nonumber\\
&-((3k+2)\alpha^{2}-(5k+3)\alpha+3k+1)p^{2}-(2\alpha^{2}-2\alpha+2)p.
\end{align}
Using Lemma 2.3 and the equitable partition $V(G)=V(K_p)\cup V(K_{n-(k+1)p-1})\cup V((kp+1)K_1)$, the largest root of $\varphi_B(x)=0$ equals $\rho_{\alpha}(G)$. Let $\beta_1=\rho_{\alpha}(G)\geq\beta_2\geq\beta_3$
be the three roots of $\varphi_B(x)=0$ and $J=\mbox{diag}(p,n-(k+1)p-1,kp+1)$. One check that
\begin{align*}
J^{\frac{1}{2}}BJ^{-\frac{1}{2}}=\left(
  \begin{array}{ccc}
  \alpha n+(1-\alpha)p-1 & (1-\alpha)p^{\frac{1}{2}}(n-(k+1)p-1)^{\frac{1}{2}} & (1-\alpha)p^{\frac{1}{2}}(kp+1)^{\frac{1}{2}}\\
  (1-\alpha)p^{\frac{1}{2}}(n-(k+1)p-1)^{\frac{1}{2}} & n+(\alpha-k-1)p-2 & 0\\
  (1-\alpha)p^{\frac{1}{2}}(kp+1)^{\frac{1}{2}} & 0 & \alpha p\\
  \end{array}
\right)
\end{align*}
is symmetric, and contains
\begin{align*}
\left(
  \begin{array}{cc}
    n+(\alpha-k-1)p-2 & 0\\
    0 & \alpha p\\
  \end{array}
\right)
\end{align*}
as a submatrix. Note that $J^{\frac{1}{2}}BJ^{-\frac{1}{2}}$ and $B$ admit the same eigenvalues. Then the Cauchy Interlacing Theorem (see Lemma 2.4) yields that $\theta_2\leq n+(\alpha-k-1)p-2<n-k-2$. Combining
this with \eqref{eq:3.1}, we deduce
\begin{align}\label{eq:3.3}
\theta=\rho_{\alpha}(G_*)>\rho_{\alpha}(K_{n-k-1})=n-k-2>\beta_2.
\end{align}
Recall that $\varphi_{B_*}(\theta)=0$. Then we possess
\begin{align}\label{eq:3.4}
\varphi_B(\theta)=\varphi_B(\theta)-\varphi_{B_*}(\theta)=(p-1)g(\theta),
\end{align}
where $g(\theta)=(k-\alpha)\theta^{2}+((\alpha^{2}-k \alpha+\alpha)n-kp-2\alpha-1)\theta-\alpha^{2}n^{2}+(2k\alpha^{2}-2k \alpha+k)pn+(k(2\alpha^{2}-2\alpha+1)+3\alpha^{2}-\alpha+1)n
-((k^{2}+2k)\alpha^{2}-(2k^{2}+3k)\alpha+k^{2}+k)p^{2}-((k^{2}+5k+2)\alpha^{2}-(2k^{2}+8k+3)\alpha+k^{2}+4k+1)p-(k^{2}+5k+4)\alpha^{2}+(2k^{2}+8k+5)\alpha-k^{2}-4k-3$. Note that the symmetry axis of $g(\theta)$
is $\theta=-\frac{(\alpha^{2}-k \alpha+\alpha)n-kp-2\alpha-1}{2(k-\alpha)}$. Then $g(\theta)$ is strictly increasing for $\theta\geq-\frac{(\alpha^{2}-k \alpha+\alpha)n-kp-2\alpha-1}{2(k-\alpha)}$. According to
\eqref{eq:3.3}, $p\geq2$, $k\geq2$, $0\leq\alpha\leq\frac{3}{4}$ and $n\geq(k+1)p+3$, we can prove $-\frac{(\alpha^{2}-k \alpha+\alpha)n-kp-2\alpha-1}{2(k-\alpha)}<n-k-2<\theta$. Thus, we have
\begin{align}\label{eq:3.5}
g(\theta)>&g(n-k-2)\nonumber\\
=&(k-\alpha)(n-k-2)^{2}+((\alpha^{2}-k \alpha+\alpha)n-kp-2\alpha-1)(n-k-2)\nonumber\\
&-\alpha^{2}n^{2}+(2k\alpha^{2}-2k \alpha+k)pn+(k(2\alpha^{2}-2\alpha+1)+3\alpha^{2}-\alpha+1)n\nonumber\\
&-((k^{2}+2k)\alpha^{2}-(2k^{2}+3k)\alpha+k^{2}+k)p^{2}\nonumber\\
&-((k^{2}+5k+2)\alpha^{2}-(2k^{2}+8k+3)\alpha+k^{2}+4k+1)p\nonumber\\
&-(k^{2}+5k+4)\alpha^{2}+(2k^{2}+8k+5)\alpha-k^{2}-4k-3\nonumber\\
=&k(1-\alpha)n^{2}+((2k\alpha^{2}-2k \alpha)p+(k+1)\alpha^{2}+(k^{2}+k-1)\alpha-2k^{2}-3k)n\nonumber\\
&-((k^{2}+2k)\alpha^{2}-(2k^{2}+3k)\alpha+k^{2}+k)p^{2}\nonumber\\
&-((k^{2}+5k+2)\alpha^{2}-(2k^{2}+8k+3)\alpha+2k+1)p\nonumber\\
&-(k^{2}+5k+4)\alpha^{2}+(k^{2}+6k+5)\alpha+k^{3}+3k^{2}+k-1\nonumber\\
\triangleq&h(n).
\end{align}

For $p\geq6$, we obtain
\begin{align}\label{eq:3.6}
-\frac{(2k\alpha^{2}-2k \alpha)p+(k+1)\alpha^{2}+(k^{2}+k-1)\alpha-2k^{2}-3k}{2k(1-\alpha)}<(k+1)p+3\leq n
\end{align}
in terms of $k\geq2$, $0\leq\alpha\leq\frac{3}{4}$ and $n\geq(k+1)p+3$. Note that the symmetry axis of $h(n)$ is $n=-\frac{(2k\alpha^{2}-2k \alpha)p+(k+1)\alpha^{2}+(k^{2}+k-1)\alpha-2k^{2}-3k}{2k(1-\alpha)}$,
which yields that $h(n)$ is strictly increasing for $n\geq-\frac{(2k\alpha^{2}-2k \alpha)p+(k+1)\alpha^{2}+(k^{2}+k-1)\alpha-2k^{2}-3k}{2k(1-\alpha)}$. Combining this with \eqref{eq:3.6}, we deduce
\begin{align*}
h(n)\geq&h((k+1)p+3)\\
=&k(1-\alpha)((k+1)p+3)^{2}\\
&+((2k\alpha^{2}-2k \alpha)p+(k+1)\alpha^{2}+(k^{2}+k-1)\alpha-2k^{2}-3k)((k+1)p+3)\\
&-((k^{2}+2k)\alpha^{2}-(2k^{2}+3k)\alpha+k^{2}+k)p^{2}\\
&-((k^{2}+5k+2)\alpha^{2}-(2k^{2}+8k+3)\alpha+2k+1)p\\
&-(k^{2}+5k+4)\alpha^{2}+(k^{2}+6k+5)\alpha+k^{3}+3k^{2}+k-1\\
=&(k^{2}\alpha^{2}-(k^{3}+2k^{2})\alpha+k^{3}+k^{2})p^{2}\\
&+((3k-1)\alpha^{2}+(k^{3}-2k^{2}-4k+2)\alpha-2k^{3}+k^{2}+k-1)p\\
&-(k^{2}+2k+1)\alpha^{2}+(4k^{2}+2)\alpha+k^{3}-3k^{2}+k-1\\
\geq&36(k^{2}\alpha^{2}-(k^{3}+2k^{2})\alpha+k^{3}+k^{2})\\
&+6((3k-1)\alpha^{2}+(k^{3}-2k^{2}-4k+2)\alpha-2k^{3}+k^{2}+k-1)\\
&-(k^{2}+2k+1)\alpha^{2}+(4k^{2}+2)\alpha+k^{3}-3k^{2}+k-1\\
=&(35k^{2}+16k-7)\alpha^{2}-(30k^{3}+80k^{2}+24k-14)\alpha+25k^{3}+39k^{2}+7k-7\\
\geq&\frac{9}{16}(35k^{2}+16k-7)-\frac{3}{4}(30k^{3}+80k^{2}+24k-14)+25k^{3}+39k^{2}+7k-7\\
=&\frac{1}{16}(40k^{3}-21k^{2}-32k-7)\\
>&0,
\end{align*}
where the last three inequalities occur from
$$
-\frac{(3k-1)\alpha^{2}+(k^{3}-2k^{2}-4k+2)\alpha-2k^{3}+k^{2}+k-1}{2(k^{2}\alpha^{2}-(k^{3}+2k^{2})\alpha+k^{3}+k^{2})}<6\leq p,
$$
$$
\frac{30k^{3}+80k^{2}+24k-14}{2(35k^{2}+16k-7)}>\frac{3}{4}\geq\alpha
$$
and $k\geq2$, respectively.

For $2\leq p\leq5$, it follows from $k\geq2$, $n\geq4k+12$ and $0\leq\alpha\leq\frac{3}{4}$ that
$$
-\frac{(2k\alpha^{2}-2k \alpha)p+(k+1)\alpha^{2}+(k^{2}+k-1)\alpha-2k^{2}-3k}{2k(1-\alpha)}<4k+12\leq n.
$$
Therefore, $h(n)$ is strictly increasing for $n\geq4k+12$. Thus, we obtain
\begin{align}\label{eq:3.7}
h(n)\geq&h(4k+12)\nonumber\\
=&k(1-\alpha)(4k+12)^{2}+((2k\alpha^{2}-2k \alpha)p+(k+1)\alpha^{2}+(k^{2}+k-1)\alpha-2k^{2}-3k)(4k+12)\nonumber\\
&-((k^{2}+2k)\alpha^{2}-(2k^{2}+3k)\alpha+k^{2}+k)p^{2}\nonumber\\
&-((k^{2}+5k+2)\alpha^{2}-(2k^{2}+8k+3)\alpha+2k+1)p\nonumber\\
&-(k^{2}+5k+4)\alpha^{2}+(k^{2}+6k+5)\alpha+k^{3}+3k^{2}+k-1\nonumber\\
=&-((k^{2}+2k)\alpha^{2}-(2k^{2}+3k)\alpha+k^{2}+k)p^{2}\nonumber\\
&+((7k^{2}+19k-2)\alpha^{2}-(6k^{2}+16k-3)\alpha-2k-1)p\nonumber\\
&+(3k^{2}+11k+8)\alpha^{2}-(12k^{3}+79k^{2}+130k+7)\alpha+9k^{3}+63k^{2}+109k-1\nonumber\\
\triangleq&l(p).
\end{align}

For $p=2$, we obtain
\begin{align*}
l(2)=&(13k^{2}+41k+4)\alpha^{2}-(12k^{3}+83k^{2}+150k+1)\alpha+9k^{3}+59k^{2}+101k-3\\
\geq&\frac{9}{16}(13k^{2}+41k+4)-\frac{3}{4}(12k^{3}+83k^{2}+150k+1)+9k^{3}+59k^{2}+101k-3\\
& \ \ \ \ \ \ \ \Big(\mbox{since} \ k\geq2 \ \mbox{and} \ 0\leq\alpha\leq\frac{3}{4}\Big)\\
=&\frac{1}{16}(65k^{2}+185k-24)\\
>&0 \ \ \ \ \ \ \ (\mbox{since} \ k\geq2).
\end{align*}

For $p=3$, we get
\begin{align*}
l(3)=&(15k^{2}+50k+2)\alpha^{2}-(12k^{3}+79k^{2}+151k-2)\alpha+9k^{3}+54k^{2}+94k-4\\
\geq&\frac{9}{16}(15k^{2}+50k+2)-\frac{3}{4}(12k^{3}+79k^{2}+151k-2)+9k^{3}+54k^{2}+94k-4\\
& \ \ \ \ \ \ \ \Big(\mbox{since} \ k\geq2 \ \mbox{and} \ 0\leq\alpha\leq\frac{3}{4}\Big)\\
=&\frac{1}{16}(51k^{2}+142k-22)\\
>&0 \ \ \ \ \ \ \ (\mbox{since} \ k\geq2).
\end{align*}

For $p=4$, we conclude
\begin{align*}
l(4)=&(15k^{2}+55k)\alpha^{2}-(12k^{3}+71k^{2}+146k-5)\alpha+9k^{3}+47k^{2}+85k-5\\
\geq&\frac{9}{16}(15k^{2}+55k)-\frac{3}{4}(12k^{3}+71k^{2}+146k-5)+9k^{3}+47k^{2}+85k-5\\
& \ \ \ \ \ \ \ \Big(\mbox{since} \ k\geq2 \ \mbox{and} \ 0\leq\alpha\leq\frac{3}{4}\Big)\\
=&\frac{1}{16}(35k^{2}+103k-20)\\
>&0 \ \ \ \ \ \ \ (\mbox{since} \ k\geq2).
\end{align*}

For $p=5$, we get
\begin{align*}
l(5)=&(13k^{2}+56k-2)\alpha^{2}-(12k^{3}+59k^{2}+135k-8)\alpha+9k^{3}+38k^{2}+74k-6\\
\geq&\frac{9}{16}(13k^{2}+56k-2)-\frac{3}{4}(12k^{3}+59k^{2}+135k-8)+9k^{3}+38k^{2}+74k-6\\
& \ \ \ \ \ \ \ \Big(\mbox{since} \ k\geq2 \ \mbox{and} \ 0\leq\alpha\leq\frac{3}{4}\Big)\\
=&\frac{1}{16}(17k^{2}+68k-18)\\
>&0 \ \ \ \ \ \ \ (\mbox{since} \ k\geq2).
\end{align*}

Thus, we have $l(p)>0$ for $2\leq p\leq5$. Together with \eqref{eq:3.7}, we admit $h(n)\geq l(p)>0$ for $2\leq p\leq5$.

From the above discussion, we conclude $h(n)>0$ for $p\geq2$. Combining this with \eqref{eq:3.4}, \eqref{eq:3.5} and $p\geq2$, we get
$$
\varphi_B(\theta)=(p-1)g(\theta)>(p-1)h(n)>0.
$$
Recall that the largest root of $\varphi_B(x)=0$ equals $\rho_{\alpha}(G)$. As $\beta_2<n-k-2<\rho_{\alpha}(G_*)=\theta$ (see \eqref{eq:3.3}), we deduce
$$
\rho_{\alpha}(G)<\rho_{\alpha}(G_*)=\rho_{\alpha}(K_1\vee(K_{n-k-2}\cup(k+1)K_1))
$$
for $p\geq2$, which contradicts $\rho_{\alpha}(G)\geq\rho_{\alpha}(K_1\vee(K_{n-k-2}\cup(k+1)K_1))$.

\noindent{\bf Case 2.} $n=(k+1)p+2$.

In this case, $G=K_p\vee(kp+2)K_1$. According to the partition $V(G)=V((kp+2)K_1)\cup V(K_p)$, the quotient matrix of $A_{\alpha}(G)$ can be expressed as
\begin{align*}
M=\left(
  \begin{array}{cc}
    \alpha p & (1-\alpha)p\\
    (1-\alpha)(kp+2) & (k\alpha+1)p+2\alpha-1\\
  \end{array}
\right).
\end{align*}
The characteristic polynomial of $M$ is
$$
\varphi_M(x)=x^{2}-((k\alpha+\alpha+1)p+2\alpha-1)x+((2k+1)\alpha-k)p^{2}+(3\alpha-2)p.
$$
Note that the partition $V(G)=V((kp+2)K_1)\cup V(K_p)$ is equitable. Then Lemma 2.3 yields that the largest root of $\varphi_M(x)=0$ equals $\rho_{\alpha}(G)$.

The symmetry axis of $\varphi_M(x)$ is $x=\frac{(k\alpha+\alpha+1)p+2\alpha-1}{2}$, and $\frac{(k\alpha+\alpha+1)p+2\alpha-1}{2}<n-k-2=(k+1)p-k$. Thus, we see that $\varphi_M(x)$ is strictly increasing for
$x\geq n-k-2$. By a simple calculation, we have
\begin{align}\label{eq:3.8}
\varphi_M(n-k-2)=&\varphi_M((k+1)p-k)\nonumber\\
=&((k+1)p-k)^{2}-((k\alpha+\alpha+1)p+2\alpha-1)((k+1)p-k)\nonumber\\
&+((2k+1)\alpha-k)p^{2}+(3\alpha-2)p\nonumber\\
=&(k^{2}-k^{2}\alpha)p^{2}+(k^{2}\alpha-k\alpha+\alpha-2k^{2}-1)p+k^{2}+2k\alpha-k\nonumber\\
\triangleq&f(p).
\end{align}
Notice that
$$
-\frac{k^{2}\alpha-k\alpha+\alpha-2k^{2}-1}{2(k^{2}-k^{2}\alpha)}<5
$$
by $k\geq2$ and $0\leq\alpha\leq\frac{3}{4}$. Then $f(p)$ is strictly increasing for $p\geq5$. Since $n=(k+1)p+2\geq4k+12$, we obtain $p\geq4+\frac{6}{k+1}$. If $k\in\{2,3,4\}$, then $p\geq6$. Hence, we get
\begin{align*}
f(p)\geq&f(6)\\
=&36(k^{2}-k^{2}\alpha)+6(k^{2}\alpha-k\alpha+\alpha-2k^{2}-1)+k^{2}+2k\alpha-k\\
=&25k^{2}-(30k^{2}+4k-6)\alpha-k-6\\
\geq&25k^{2}-\frac{3}{4}(30k^{2}+4k-6)-k-6\\
=&\frac{1}{2}(5k^{2}-8k-3)\\
>&0.
\end{align*}
If $k\geq5$, then $p\geq5$. Hence, we possess
\begin{align*}
f(p)\geq&f(5)\\
=&25(k^{2}-k^{2}\alpha)+5(k^{2}\alpha-k\alpha+\alpha-2k^{2}-1)+k^{2}+2k\alpha-k\\
=&16k^{2}-(20k^{2}+3k-5)\alpha-k-5\\
\geq&16k^{2}-\frac{3}{4}(20k^{2}+3k-5)-k-5\\
=&\frac{1}{4}(4k^{2}-13k-5)\\
>&0.
\end{align*}
From the above discussion, we infer $f(p)>0$. Together with \eqref{eq:3.8}, we conclude $\varphi_M(n-k-2)>0$. Recall that $\varphi_M(x)$ is strictly increasing for $x\geq n-k-2$ and the largest root of $\varphi_M(x)=0$
equals $\rho_{\alpha}(G)$. This yields that $\rho_{\alpha}(G)<n-k-2$. Then by \eqref{eq:3.1}, we admit
$$
\rho_{\alpha}(G)<\rho_{\alpha}(G_*)=\rho_{\alpha}(K_1\vee(K_{n-k-2}\cup(k+1)K_1)),
$$
which contradicts $\rho_{\alpha}(G)\geq\rho_{\alpha}(K_1\vee(K_{n-k-2}\cup(k+1)K_1))$.

\noindent{\bf Case 3.} $n=(k+1)p+1$.

In this case, $G=K_p\vee(kp+1)K_1$. In light of the partition $V(G)=V((kp+1)K_1)\cup V(K_p)$, the quotient matrix of $A_{\alpha}(G)$ is
\begin{align*}
R=\left(
  \begin{array}{cc}
    \alpha p & (1-\alpha)p\\
    (1-\alpha)(kp+1) & (k\alpha+1)p+\alpha-1\\
  \end{array}
\right).
\end{align*}
The characteristic polynomial of $R$ is
$$
\varphi_R(x)=x^{2}-((k\alpha+\alpha+1)p+\alpha-1)x+((2k+1)\alpha-k)p^{2}+(\alpha-1)p.
$$
Since the partition $V(G)=V((kp+1)K_1)\cup V(K_p)$ is equitable, it follows from Lemma 2.3 that the largest root of $\varphi_R(x)=0$ equals $\rho_{\alpha}(G)$.

Since $n=(k+1)p+1\geq4k+12$, we possess $p\geq4+\frac{7}{k+1}$. The symmetry axis of $\varphi_R(x)$ is $x=\frac{(k\alpha+\alpha+1)p+\alpha-1}{2}$, and $\frac{(k\alpha+\alpha+1)p+\alpha-1}{2}<n-k-2=(k+1)p-k-1$. This
implies that $\varphi_R(x)$ is strictly increasing for $x\geq n-k-2$. A simple computation yields that
\begin{align}\label{eq:3.9}
\varphi_R(n-k-2)=&\varphi_R((k+1)p-k-1)\nonumber\\
=&((k+1)p-k-1)^{2}-((k\alpha+\alpha+1)p+\alpha-1)((k+1)p-k-1)\nonumber\\
&+((2k+1)\alpha-k)p^{2}+(\alpha-1)p\nonumber\\
=&(k^{2}-k^{2}\alpha)p^{2}+(k^{2}\alpha+k\alpha+\alpha-2k^{2}-2k-1)p+k^{2}+k+k\alpha+\alpha\nonumber\\
\triangleq&w(p).
\end{align}
Note that
$$
-\frac{k^{2}\alpha+k\alpha+\alpha-2k^{2}-2k-1}{2(k^{2}-k^{2}\alpha)}<5
$$
by $k\geq2$ and $0\leq\alpha\leq\frac{3}{4}$. Then $w(p)$ is strictly increasing for $p\geq5$. If $k=2$, then $p\geq7$. Thus, we possess
\begin{align*}
w(p)\geq&w(7)\\
=&49(k^{2}-k^{2}\alpha)+7(k^{2}\alpha+k\alpha+\alpha-2k^{2}-2k-1)+k^{2}+k+k\alpha+\alpha\\
=&36k^{2}-13k-7-(42k^{2}-8k-8)\alpha\\
\geq&36k^{2}-13k-7-\frac{3}{4}(42k^{2}-8k-8)\\
=&\frac{1}{2}(9k^{2}-14k-2)\\
>&0.
\end{align*}
If $k\in\{3,4,5\}$, then $p\geq6$. Thus, we conclude
\begin{align*}
w(p)\geq&w(6)\\
=&36(k^{2}-k^{2}\alpha)+6(k^{2}\alpha+k\alpha+\alpha-2k^{2}-2k-1)+k^{2}+k+k\alpha+\alpha\\
=&25k^{2}-11k-6-(30k^{2}-7k-7)\alpha\\
\geq&25k^{2}-11k-6-\frac{3}{4}(30k^{2}-7k-7)\\
=&\frac{1}{4}(10k^{2}-23k-3)\\
>&0.
\end{align*}
If $k\geq6$, then $p\geq5$. Thus, we have
\begin{align*}
w(p)\geq&w(5)\\
=&25(k^{2}-k^{2}\alpha)+5(k^{2}\alpha+k\alpha+\alpha-2k^{2}-2k-1)+k^{2}+k+k\alpha+\alpha\\
=&16k^{2}-9k-5-(20k^{2}-6k-6)\alpha\\
\geq&16k^{2}-9k-5-\frac{3}{4}(20k^{2}-6k-6)\\
=&\frac{1}{2}(2k^{2}-9k-1)\\
>&0.
\end{align*}
From the above discussion, we obtain $w(p)>0$. Combining this with \eqref{eq:3.9}, we possess $\varphi_R(n-k-2)>0$. Recall that $\varphi_R(x)$ is strictly increasing for $x\geq n-k-2$ and the largest root of
$\varphi_R(x)=0$ equals $\rho_{\alpha}(G)$. Hence, we deduce $\rho_{\alpha}(G)<n-k-2$. Then by \eqref{eq:3.1}, we conclude
$$
\rho_{\alpha}(G)<\rho_{\alpha}(G_*)=\rho_{\alpha}(K_1\vee(K_{n-k-2}\cup(k+1)K_1)),
$$
which contradicts $\rho_{\alpha}(G)\geq\rho_{\alpha}(K_1\vee(K_{n-k-2}\cup(k+1)K_1))$. Theorem 1.2 is proved. \hfill $\Box$

\section*{Declaration of competing interest}

The author declares that he has no known competing financial interests or personal relationships that could have appeared to influence the work reported in this paper.

\section*{Data availability}

No data was used for the research described in the article.

\section*{AI Contribution}

No use of AI tools were made in the preparation of this article.



\section*{Acknowledgments}

This work was supported by the Natural Science Foundation of Jiangsu Province (Grant No. BK20241949).

\end{document}